\documentclass{gtart}
\def\ifplaintex{\expandafter\ifx\csname documentclass\endcsname\relax}

\def\gtp{{\mathsurround=0pt\it $\cal G\mskip-2mu$eometry \&\ 
$\cal T\!\!$opology $\cal P\!$ublications}}  % GT publications

\def\recd{{\small Received:\qua\receiveddate\ifx\reviseddate\relax
\else\qquad Revised:\qua\reviseddate\fi\par}} 

\def\lognumber#1{\def\thelognumber{#1}}
\def\volumenumber#1{\def\thevolumenumber{#1}}
\def\volumeyear#1{\def\thevolumeyear{#1}}
\def\papernumber#1{\def\thepapernumber{#1}}
\def\pagenumbers#1#2{\def\startpage{#1}\def\finishpage{#2}}
\def\published#1{\def\publishdate{#1}}

\def\received#1{\def\receiveddate{#1}}

\def\accepted#1{\def\accepteddate{#1}}

\long\def\asciiabstract#1{\long\def\theasciiabstract{#1}}

\let\\\par\let\thelognumber\relax\let\thevolumenumber\relax
\let\thepapernumber\relax\let\thevolumeyear\relax\let\startpage\relax
\let\finishpage\relax\let\publishdate\relax\let\receiveddate\relax
\let\reviseddate\relax\let\accepteddate\relax\let\theasciititle\relax
\let\theasciiauthors\relax
\let\theasciiabstract\relax

\let\theasciiemail\relax

\ifplaintex
\font\logobig=cmssbx10 scaled 3836
\font\logomed=cmssbx10 scaled 2557
\else
\font\logobig=cmssbx10 scaled 4200
\font\logomed=cmssbx10 scaled 2800
\fi

\long\def\makeagttitle{   %%% start of definition of \makeagttitle
\count0=\startpage
\agt\hfill      %   Journal title (top left) 
\hbox to 45truept{\vbox to 0pt{\vglue -13truept{\logomed A\kern -.37em{\logobig 
T}\kern -.38em G}\vss}\hss}
\break
{\small Volume \thevolumenumber\ (\thevolumeyear)
\startpage--\finishpage\nl
Published: \publishdate}

\vglue .25truein

{\parskip=0pt\leftskip 0pt plus
1fil\def\\{\par\smallskip}{\Large\bf\thetitle}\par\medskip} \vglue
0.05truein

{\parskip=0pt\leftskip 0pt plus 1fil\def\\{\par}{\sc\theauthors}
\par\medskip}%
 
\vglue 0.03truein 

{\small\leftskip 25truept\rightskip 25truept{\bf Abstract}\stdspace\theabstract

{\bf AMS Classification}\stdspace\theprimaryclass
\ifx\thesecondaryclass\relax\else; \thesecondaryclass\fi\par
{\bf Keywords}\stdspace \thekeywords\par}\vglue 7truept

}   %%%% end of definition of \makeagttitle

\ifplaintex
\font\phead=cmsl9 scaled 950
\font\pnum=cmbx10 scaled 913
\font\pfoot=cmsl9 scaled 950
\headline{\vbox to 0pt{\vskip -4.5mm\line{\small\phead\ifnum
\count0=\startpage ISSN 1472-2739 (on-line) 1472-2747 (printed)
\hfill {\pnum\folio}\else\ifodd\count0\def\\{ }% 
\ifx\theshorttitle\relax\thetitle\else\theshorttitle\fi\hfill{\pnum\folio}
\else\def\\{ and }{\pnum\folio}\hfill\ifx\theshortauthors\relax\theauthors
\else\theshortauthors\fi\fi\fi}\vss}}
\footline{\vbox to 0pt{\vglue 0mm\line{\small\pfoot\ifnum\count0=\startpage
\copyright\ \gtp\hfill\else
\agt, Volume \thevolumenumber\ (\thevolumeyear)\hfill\fi}\vss}}
\else
\font=cmsl9 scaled 1050
\font\lnum=cmbx10 
\font=cmsl9 scaled 1050
\makeatletter
\def\@oddhead{{\small\lhead\ifnum\count0=\startpage ISSN 1472-2739 
(on-line) 1472-2747 (printed)\hfill {\lnum\number\count0}\else\ifodd\count0
\def\\{ }\ifx\theshorttitle\relax \thetitle \else\theshorttitle\fi\hfill
{\lnum\number\count0}\else\def\\{ and }{\lnum\number\count0}
\hfill\ifx\theshortauthors\relax 
\theauthors\else\theshortauthors\fi\fi\fi}}\def\@evenhead{\@oddhead}
\def\@oddfoot{\small\lfoot\ifnum\count0=\startpage\copyright\ \gtp\hfill\else
\agt, Volume \thevolumenumber\ (\thevolumeyear)\hfill\fi}
\def\@evenfoot{\@oddfoot}
\makeatother
\fi
\let\maketitlepage\makeagttitle
\let\makeshorttitle\maketitlepage
\let\maketitle\maketitlepage

\newwrite\gtoutfile
\long\gdef\makeheadfile{  %%% start of definition of \makeheadfile
{\def\\{, }\def\s{ }
\immediate\openout\gtoutfile head.xxx
\immediate\write\gtoutfile{To: math@arxiv.org}
\immediate\write\gtoutfile{Subject: put OR rep NNNNN:ppppp}
\immediate\write\gtoutfile{--text follows this line--}
\immediate\write\gtoutfile{Proxy-for: \ifx\theasciiauthors\relax
\theauthors\else\theasciiauthors\fi\s<\ifx\theasciiemail\relax\theemail\else\theasciiemail\fi>}
\immediate\write\gtoutfile{\noexpand\\}
\immediate\write\gtoutfile{Authors: \ifx\theasciiauthors\relax
\theauthors\else\theasciiauthors\fi}
{\def\\{ }\immediate\write\gtoutfile{Title: \ifx\theasciititle\relax
\thetitle\else\theasciititle\fi}}
\immediate\write\gtoutfile{Subj-class: GT or SG, GR etc}
\immediate\write\gtoutfile{MSC-class: \theprimaryclass\ifx\thesecondaryclass\relax\else, \thesecondaryclass\fi}
\immediate\write\gtoutfile{Journal-ref: Algebr. Geom. Topol. \thevolumenumber\s
(\thevolumeyear) \startpage-\finishpage}
\immediate\write\gtoutfile{Comments: Published by Algebraic and
Geometric Topology at}
\immediate\write\gtoutfile{\s\s\s  http://www.maths.warwick.ac.uk/agt/AGTVol\thevolumenumber/agt-\thevolumenumber-\thepapernumber.abs.html}
\immediate\write\gtoutfile{\noexpand\\}
\immediate\write\gtoutfile{}
\ifx\theasciiabstract\relax
\immediate\write\gtoutfile{\theabstract}\else
\immediate\write\gtoutfile{\theasciiabstract}\fi
\immediate\write\gtoutfile{}
\immediate\write\gtoutfile{\noexpand\\}
\immediate\write\gtoutfile{}
\immediate\closeout\gtoutfile}}  %%% end of definition of \makeheadfile

\def\maketitlepage{\makeagttitle\makeheadfile}
\let\makeshorttitle\maketitlepage
\let\maketitle\maketitlepage

\lognumber{47}
\volumenumber{2}
\volumeyear{2002}
\papernumber{47}
\pagenumbers{1197}{1204}
\received{5 August 2002}
\accepted{17 December 2002}
\published{27 December 2002}

\usepackage{amssymb,amsmath}
\usepackage{pb-diagram,lamsarrow,pb-lams}
\usepackage[rokicki,hideboxes]{boxedeps}

\def\lbl#1{\label{#1}}
\theoremstyle{plain}

\newtheorem{theorem}{Theorem}
\newtheorem{proposition}{Proposition}[section]
\newtheorem{lemma}[proposition]{Lemma}
\newtheorem{corollary}[proposition]{Corollary}

\newtheorem{conjecture}{Conjecture}

\theoremstyle{definition}

\newcommand{\psdraw}[2]{\centerline{\BoxedEPSF{#1.eps scaled
#2}}}
\newlength{\standardunitlength}
\renewcommand{\ker}{\operatorname{Ker}}
\newcommand{\im}{\operatorname{Im}}

\def\ov{\overset}
\def\ak{\mathcal A_k^t (H)}

\def\SS{\Sigma}

\def\a{\alpha}
\def\b{\beta}

\def\th{\theta}

\def\t{\tau}
\newcommand{\s}{\sigma}
\def\g{\gamma}

\def\Z{\mathbb Z}

\def\D{\mathsf D}

\def\A{\mathcal A}

\def\G{\mathcal G}
\def\Q{\mathbb Q}

\def\Sg{\Sigma_{g,1}}

\def\iso{\cong}

\def\sub{\subseteq}

\def\Gg{\Gamma_{g,1}}

\def\Hg{\mathcal H_g}

\def\dk{\D_k (H)}
\def\dki{\D'_k (H)}
\def\ll#1{L_{#1}(H)}
\def\fw#1{\mathcal F^w_{#1}(\Gg )}
\def\gw#1{\mathcal G^w_{#1}(\Gg )} 
\def\fy#1{\mathcal F^Y_{#1}(\Hg )}
\def\gy#1{\mathcal G^Y_{#1}(\Hg )} 
\def\fwh#1{\mathcal F^w_{#1}(\Hg )}
\def\gwh#1{\mathcal G^w_{#1}(\Hg )}

\begin{document}
\title
[Addendum and correction to: Homology cylinders]{Addendum and correction to:\\Homology
cylinders:\\
an enlargement of the mapping
class group}
\author{Jerome Levine}
\address{Department of Mathematics, Brandeis University\\Waltham,
 MA 02454-9110, USA}        
\email{levine@brandeis.edu}
\url{http://people.brandeis.edu/\char'176levine/} 

\begin{abstract}
In a previous paper  \cite{Le}, a group $\Hg$ of
{\em
homology cylinders} over the oriented surface of genus $g$ is defined. A
filtration of $\Hg$ is defined, using the Goussarov-Habiro notion of
finite-type. It is erroneously claimed that this filtration essentially
coincides with the relative weight filtration. The present note corrects
this error and studies the actual relation between the two filtrations.
\end{abstract}

\asciiabstract{In a previous paper [Homology cylinders: an enlargement
of the mapping class group, Algebr. Geom. Topol. 1 (2001) 243--270,
arXiv:math.GT/0010247], a group H_g of homology cylinders over the
oriented surface of genus g is defined. A filtration of H_g is
defined, using the Goussarov-Habiro notion of finite-type. It is
erroneously claimed that this filtration essentially coincides with
the relative weight filtration. The present note corrects this error
and studies the actual relation between the two filtrations.}


\maketitle

\section{Introduction}
In \cite{Le} we consider a group $\Hg$ consisting of homology bordism
classes of {\em homology cylinders}, where a homology cylinder is
defined as a homology bordism between two copies of $\Sg$, the once
punctured oriented surface of genus $g$. This bordism is equipped with
an
explicit identification of each end with $\Sg$---see \cite{Le} for more
details. In particular there is a canonical injection of the mapping
class group $\Gg$ into $\Hg$.

Two filtrations of $\Hg$ are considered in the first part of the paper.
The first
is the relative weight filtration $\fwh{k}$, the obvious extension
of the relative weight filtration $\fw{k}$ of $\Gg$ considered by
Johnson, Morita and others. The canonical injection $J_k
:\gw{k}\to\dk$, where $\gw{k}=\fw{k}/\fw{k+1}$ is the associated
graded group,
lifts, in a natural way, to $\gwh{k}$ and there becomes an {\em
isomorphism} $J_k^H :\gwh{k}\overset{\iso}\longrightarrow\dk$. $\dk$
is
the kernel of the bracket map $\b_k :H\otimes\ll{k+1}\to\ll{k+2}$, where
$\ll{k}$ is the degree $k$ component of the free Lie algebra on $H=H_1
(\Sg )$.

The second filtration $\fy{k}$ is defined using the Goussarov-Habiro
theory of finite-type invariants of $3$-manifolds. It is shown in
\cite{Le} that $\fy{k}\sub\fwh{k}$ and so we have induced homomorphisms
$\gy{k}\to\gwh{k}\iso\dk$. To study these homomorphisms we use results
{\em
announced} by Habiro. Habiro considers an abelian group $\A_k (H)$
defined by unitrivalent graphs with $k$ trivalent vertices and with
univalent vertices labelled by elements of $H$, subject to
anti-symmetry, the IHX relation and linearity of labels (see \cite{Le}
for a more complete description).
We then consider the quotient $\ak$ in which only trees are allowed.
Using results of Habiro it is proved in \cite{Le} that there is a
well-defined {\em epimorphism} $\th_k :\ak\to\gy{k}$. Furthermore there
is a combinatorially defined homomorphism $\eta_k :\ak\to\dk$ (which can
be defined for {\em
any} abelian group $H$) which coincides with the
composition:
\begin{equation}\lbl{eq.comp}
\begin{diagram}
\node{\ak}\arrow{e,t,A}{\th_k}\node{\gy{k}}\arrow{e}\node{\gwh{k}}\arrow{e,tb}{J_k
(H)}{\iso}\node{\dk}
\end{diagram}
\end{equation}
Note that this is different from the map called $\eta_k$ in \cite{Le}.

It is erroneously claimed in Proposition 2.2 of \cite{Le} that $\eta_k$
is an isomorphism for $k>1$. But in fact this is FALSE. Thus the
implications that all the maps in diagram \eqref{eq.comp} are
isomorphisms for $k>1$ is false. (For $k=1$ the result
$\gy{1}\iso\gwh{1}\oplus V$, where the projection $\gy{1}\to V$ is
defined by Birman-Craggs homomorphisms, is still true.) 

It is the aim of this note to correct this error and, in
particular, study the homomorphism $\eta_k$. In fact it is known, and
will be reproved below, that $\eta_k$ induces an isomorphism
$\ak\otimes\Q\iso\dk\otimes\Q$. Thus the maps in diagram \eqref{eq.comp}
are isomorphisms for $k>1$ after tensoring with $\Q$. To handle the more general
case it will be natural to introduce a variation on the notion of Lie
algebra by replacing the axiom $[x,x]=0$ with the weaker anti-symmetry
axiom $[x,y]+[y,x]=0$ and investigate the corresponding free
objects. This variation does not seem to have been studied before, even
though it arises naturally from the study of oriented graphs.

This work was partially supported by an NSF grant
           and by an Israel-US BSF grant.

\section{A different notion of Lie algebra}

We want to discuss the map $\eta_k :\A_k^t (H)\to\D_k (H)$, for $k>1$.
For this purpose it will be more appropriate to
consider a replacement for the free Lie algebra $L(H)$. Let us define a
{\em quasi-Lie algebra} by replacing the axiom $[x,x]=0$ with the axiom
$[x,y]+[y,x]=0$, for any $x,y$. Thus we only can conclude $2[x,x]=0$,
and so if $L$ is quasi-Lie algebra then $L\otimes\Z [1/2]$ is a Lie
algebra. We can now define the free quasi-Lie algebra $L'(H)$ over a
free abelian group $H$ in the obvious way (using the free {\em magma}
over $H$, for example---see \cite{R}). There is an obvious map $\g
:L'(H)\to L(H)$, which is a map of quasi-Lie algebras. Let $\g_k :L'_k
(H)\to L_k (H)$ be the degree $k$ component.

\begin{lemma}\lbl{lem.quasi}
\begin{enumerate}
\item If $k$ is odd then $\g_k$ is an isomorphism.
\item If $k=2l$, then there is an exact sequence of additive
homomorphisms
$$
\begin{diagram}
\dgARROWLENGTH=1.2em
\node{L_l (H)/2L_l (H)}\arrow{e}\node{L'_k
(H)}\arrow{e,t}{\g_k}\node{L_k (H)\to 0}
\end{diagram}
$$
\end{enumerate}
\end{lemma}

\begin{proof} Clearly $\g_k$ is onto. Furthermore the kernel $K_k (H)$
of
$\g_k$ is generated additively by all brackets which contain a
sub-bracket of the form $[\a ,\a ]$ for some $\a\in L'(H)$. In fact such
a bracket will be zero in $L'(H)$ unless it is exactly of the form $[\a
,\a ]$. In other words for any $\a ,\eta\in L'(H)$
$$[[\a ,\a ],\eta ]=0=[\eta ,[\a ,\a ]] $$
This follows directly from the Jacobi relation and anti-symmetry. 

Thus we can define a map $L'(H)\to L'(H)$ by $\a\longmapsto [\a ,\a
]$---it is an additive homomorphism by anti-symmetry--- and the image of
this map is exactly the kernel of $\g$. Note that this map vanishes on
$2L'(H)$ and on any element of the form $\a =[\eta ,\eta ]$. The
assertions
of the lemma follow.
\end{proof}

\begin{conjecture} It is easy to see that $L_l(H)/2L_l(H)\to
L'_{2l}(H)$
is a monomorphism for $l=1$ and it is reasonable to conjecture that this
is true for all $l$.
\end{conjecture}

Analogous to $\b_k$ we can define a homomorphism $\b'_k :H\otimes
L'_{k+1}
(H)\to L'_{k+2}(H)$ by $\b'_k (h\otimes\a )=[h,\a ]$. We see that
$\b'_k$ is
onto by the Jacobi identity and anti-symmetry and denote the kernel by
$\D'_k (H)$. If we apply the snake lemma to the diagram:
$$
\dgARROWLENGTH=1.1em
\begin{diagram}
%\node{0}\arrow{e}
\node{0\to\dki}\arrow{e}\arrow{s}\node{H\otimes
L'_{k+1}(H)}\arrow{e,t}{\b'_k}\arrow{s,l}{1\otimes\g_{k+1}}\node{L'_{k+2}(H)\to 0}
%\arrow{e}
\arrow{s,l}{\g_{k+2}}%\node{0}
\\
%\node{0}\arrow{e}
\node{0\to\dk}\arrow{e}\node{H\otimes
L_{k+1}(H)}\arrow{e,t}{\b_k}\node{L_{k+2}(H)\to 0}
%\arrow{e}\node{0}
\end{diagram}
$$
 we
conclude:

\begin{corollary}\lbl{cor.dd}
The canonical map $\D'_k (H)\to\D_k (H)$ fits into the following exact
sequences, depending on whether $k$ is odd or even.
$$ 0\to\D'_{2l}(H)\to\D_{2l}(H)\to K_{2l+2}(H)\to 0 $$
$$ 0\to H\otimes K_{2l}(H)\to\D'_{2l-1}(H)\to\D_{2l-1}(H)\to 0$$
\end{corollary}

\section{$\ak$ and Lie algebras}

We will refer to a univalent vertex of a tree as a {\em leaf}, except
when the tree is {\em rooted}, i.e. one of the univalent vertices is
designated a root. In that case only the remaining univalent vertices
will be referred to as leaves. 

We can graphically interpret $L_k'(H)$ as the abelian group
generated by rooted binary planar trees with $k$ leaves, whose leaves
are labelled by
elements of $H$ modulo the anti-symmetry and IHX relations and linearity
of labels. These
relations correspond exactly to the axioms for a quasi-Lie algebra. The
correspondence is described in \cite{R}, for the case of a free magma.
Similarly we can interpret $H\otimes L'_k (H)$ as rooted binary planar
trees with $k$ leaves whose leaves {\em and root} are labelled by
elements of $H$, modulo anti-symmetry, IHX and label linearity. See
Figure
\ref{fig.tree}. 

\begin{figure}[ht!]
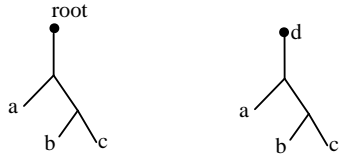

\centering
\psdraw{tree}{700}
\bigskip
\caption{The left-hand tree corresponds to $[a,[b,c]]$ in $L'_3 (H)$.
The right-hand tree corresponds to $d\otimes [a,[b,c]]$ in
$H\otimes
L'_3 (H)$.}\lbl{fig.tree}
\end{figure}

We now define
a map $\eta'_k :\ak\to H\otimes L'_{k+1}(H)$, in the same way that
$\eta_k$
was
defined, by sending each labelled binary planar tree to the sum of the
rooted labelled binary planar trees, one for each leaf, obtained by
designating that leaf as the (labelled) root. We want to show that $\im
\b'_k =\dki$, i.e. that the following sequence is exact.
$$ \ak\ov{\eta'_k}\longrightarrow H\otimes
L'_{k+1}(H)\ov{\b'_k}\longrightarrow L'_{k+2}(H)\to 0 $$
We first prove:
\begin{lemma}\lbl{lem.root}  $\im\eta'_k\sub\D'_k (H)$.
\end{lemma}
\begin{proof}
Let $T$ be a labelled planar binary tree, representing $t\in\ak$. Then
$\b'_k\circ\eta'_k (t)\in L'_{k+2}(H)$ is represented by a sum
$\sum_l
T_l$\ , over all leaves $l$ of $T$, where $T_l$ is the rooted tree
obtained by adjoining to the edge of $T$ containing $l$ a rooted edge as
in Figure \ref{fig.root}. 
We need to show that this sum represents $0$.

\begin{figure}[ht!]
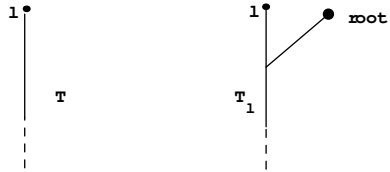

\psdraw{root}{500}
\caption{Define a rooted tree from a tree and one of its
leaves}\lbl{fig.root}
\end{figure}

Consider now the sum $\sum_{(v,e)}T_{v,e}$\ , over all pairs $(v,e)$,
where $v$ is a vertex of $T$ and $e$ an edge containing $v$. $T_{v,e}$
is the rooted tree obtained by adjoining to $e$, near $v$, a rooted edge
as in Figure \ref{fig.root1}. 
The terms of this sum for univalent
vertices $v$ is clearly just $\sum_l T_l$. The remaining terms
correspond to the internal vertices and, for each internal vertex, there
are three terms, whose
sum will vanish by the IHX relation.
Thus
it suffices to prove that $\sum_{(v,e)}T_{v,e}$ represents $0$. But
this is clear since, for each edge $e$, with vertices $v' ,v''$, we have
$T_{v',e}=-T_{v'',e}$, by the anti-symmetry relation.
\end{proof} 

\begin{figure}[ht!]
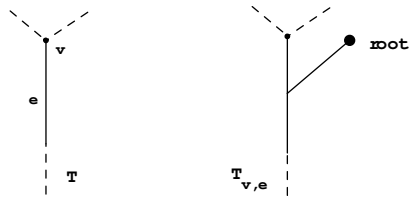

\psdraw{root1}{500}
\caption{Define a rooted tree from a tree and an edge-vertex
pair}\lbl{fig.root1}
\end{figure}

\begin{theorem}\lbl{th.tree}
$\eta '_k :\ak\to\dki$ is a split surjection. $\ker\eta '_k$ is the
torsion subgroup
of
$\ak$, if $k$ is even. It is the odd torsion subgroup if $k$ is odd. In
either case 
$$(k+2)\ker\eta '_k =0.$$
\end{theorem}

\begin{corollary} $\ak\otimes\Q\iso\D'_k (H)\otimes\Q\iso\D_k
(H)\otimes\Q $
\end{corollary}
\begin{conjecture}
It is reasonable to conjecture that $\eta'_k$ is an isomorphism.
\end{conjecture}

\begin{proof}[Proof of Theorem \ref{th.tree}] We will need some
auxiliary maps. First we
define
$$\rho_k :H\otimes L'_{k+1}(H)\to\ak .$$
Let $\a$ be a generator of $H\otimes L'_{k+1}(H)$ represented by a
rooted tree with labels on
all leaves and root. Define $\rho_k (\a )$ to be the same labelled tree
obtained by just 
forgetting which vertex is the root. This obviously preserves the
anti-symmetry and IHX
relations and label linearity, and so gives
a well-defined additive homomorphism. The important property to observe
is
$$\rho_k\circ\eta'_k =\text{multiplication by }k+2 .$$
This shows that $(k+2)\ker\eta'_k =0$.

When $k$ is even, $\D'_k (H)$ is torsion-free, by Corollary
\ref{cor.dd}. This shows
that $\ker\eta'_k$ is
the torsion subgroup of $\ak$. If $k$ is odd, then Corollary
\ref{cor.dd} shows that all the
torsion in $\D'_k (H)$ is of order 2, since $L(H)$ is torsion-free and $K_{2l}(H)$ is 2-torsion by Lemma \ref{lem.quasi}. Since $k+2$ is odd, we conclude
that $\ker\eta'_k$ is the
odd torsion subgroup of $\ak$. From this it follows that $\eta'_k$
splits.

It remains only to show that $\eta'_k$ is onto. For this we will
construct
a map
$$\t_k :L'_{k+2}(H)\to H\otimes L'_{k+1}(H)/\eta'_k (\ak ).$$
Consider a generator $\a$ of $L'_{k+2}(H)$ represented by a rooted tree
$T$ with labelled
leaves as in Figure \ref{fig.beta}. Here $v$ is the trivalent vertex
adjacent to the
root and $A, B$ the
two subtrees (rooted, with labelled leaves) with $v$ as their common
root. We then form a labelled 
tree $T'$ from $A$ and $B$ by eliminating the root of $T$ and making $v$
the midpoint of an
edge connecting $A$ to $B$---see Figure \ref{fig.beta} 

\begin{figure}[ht]
\centering
\psdraw{beta}{500}
\bigskip
\nocolon\caption{}\lbl{fig.beta}
\end{figure}

Now for each leaf $w$ of $T'$
we can create a
rooted tree $T_w$ by making $w$ the root. Then $T_w$ represents an
element of $H\otimes
L'(H)$. Recall that we defined $\eta'_k (T')$ to be the sum $\SS_w T_w$
over all leaves of $T'$. We now define $\t_k (\a )$ to
be the class represented by the sum $\SS_w T_w$
 over all
leaves  $w$ of $A$. We need to check that this is
well-defined modulo
$\eta'_k (\ak )$. 

If we consider an anti-symmetry relation in $T$ at a trivalent vertex in
$A$ or $B$ then the image is clearly an anti-symmetry relation in every
$T_w$. The anti-symmetry relation at the vertex $v$ is easily seen to
map to precisely $\eta'_k (T')$. 

Now consider an IHX relation at an internal edge $e$ of $T$. If $e$ is
an internal edge of $A$ or $B$ then it induces an IHX relation in every
$T_w$. Suppose, on the other hand that $e$ contains $v$. If the other
vertex of $e$ is in $A$, for example, then we can represent the IHX
relation  as in
Figure \ref{fig.ihx}.

\begin{figure}[ht!]
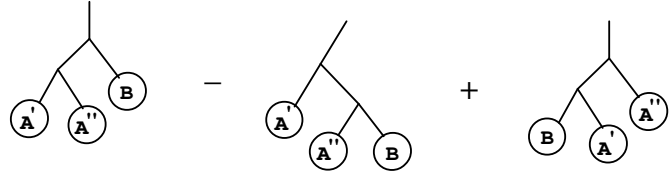

\centering
\psdraw{ihx}{700}
\bigskip
\caption{Graphical representation of the IHX relation}\lbl{fig.ihx}
\end{figure}

 Here we have split $A$ into two subtree pieces $A'$ and $A''$.
The image of this IHX relation is pictured in Figure \ref{fig.ihxim},
where we take the sum over all leaves $w$ in each subtree.

\begin{figure}[ht!]
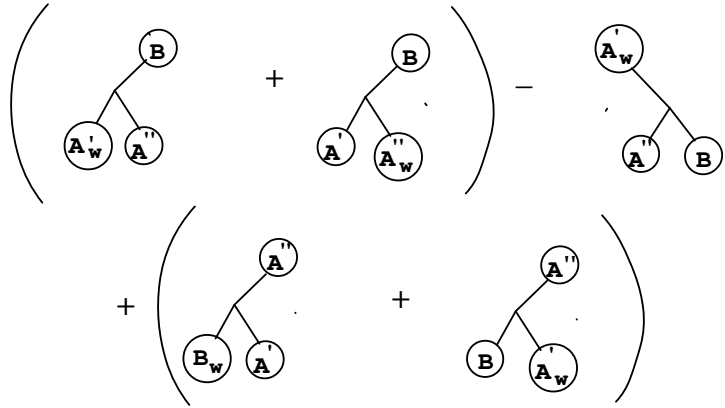

\centering
\psdraw{ihxim}{700}
\bigskip
\caption{The image of the IHX relation in Figure
\ref{fig.ihx}}\lbl{fig.ihxim}
\end{figure}

 Then we can see
that the first and third terms cancel while the second, fourth and fifth
terms
add up to exactly $\eta'_k (T')$.

Now it is easy to see that the composition 
$$ \t_k\circ\b'_k :H\otimes
L_{k+1}'(H)\to H\otimes
L_{k+1}'(H)/\eta'_k
(\ak )
 $$
is just the canonical projection. From this it follows immediately that
$\ker\b'_k =\im\eta'_k$.
\end{proof}

\section{Relation between $\gy{k}$ and $\gwh{k}$}

Finally we can draw some conclusions about the natural map $ \G_k^Y (\Hg
)\to \G_k^w (\Hg )$.

\begin{corollary}
\begin{enumerate}
\item For all $k$
$$\ak\otimes\Q\iso \G_k^Y (\Hg )\otimes\Q\iso \G_k^w (\Hg
)\otimes\Q$$
\item For $k=1$ we have $\G_1^Y (\Hg )\iso V\oplus \G_1^w (\Hg )$.
\item If $k$ is even, there is an exact sequence
$$ \G_k^Y (\Hg )\to \G_k^w (\Hg )\to K_{k+2}(H)\to 0$$
\item If $k>1$ is odd, then $ \G_k^Y (\Hg )\to \G_k^w (\Hg )$ is onto
and
there is an exact sequence
$$H\otimes K_{k+1}(H)\to\G_k^Y (\Hg )\otimes\Z_{(2)}\to \G_k^w (\Hg )
\otimes\Z_{(2)}\to 0$$
where $\Z_{(2)}$ is the ring of fractions with odd denominator.
\end{enumerate}
\end{corollary}

\begin{conjecture}
Taking account of the various conjectures mentioned above we can
conjecture that the precise relationship between $\G_k^w (\Hg )$ and
$\G_k^Y (\Hg )$, for $k>1$ is given by the following exact sequences:
\begin{gather*}
0\to \G_{2l}^Y (\Hg )\to \G_{2l}^w (\Hg )\to L_{l+1}(H)/2L_{l+1}(H)\to
0\\
H\otimes L_{l}(H)/2L_{l}(H)\to\G_{2l-1}^Y (\Hg )\to \G_{2l-1}^w (\Hg
)\to 0\quad (l>1)\end{gather*}
\end{conjecture}

\Addresses\recd
\end{document}